\documentclass[11pt]{article}
\usepackage{amssymb}
\usepackage{amsmath}
\numberwithin{equation}{section}
\usepackage{amsthm}
\usepackage{amsfonts}
\usepackage{mathrsfs}
\usepackage{makeidx}
\usepackage{multicol}

\def\qed{\hfill $\Box$}
\newcommand{\overbar}[1]{\mkern 1mu\overline{\mkern-1mu#1\mkern-1mu}\mkern 1mu}

\begin{document}
\title{\bf On the Number of Zeros and Poles of Dirichlet Series}
\author{Bao Qin Li }
\date{}
\maketitle

\bigskip
\noindent
\begin{abstract}
\noindent This paper investigates lower bounds on the number of zeros and poles of a general Dirichlet series in a disk of radius $r$ and gives, as a consequence, an affirmative answer to an open problem of Bombieri and Perelli on the bound. Applications will also be given to Picard type theorems, global estimates on the symmetric difference of zeros, and uniqueness problems for Dirichlet series.
\end{abstract}

\pagenumbering{arabic}

\begin{figure}[b]
\rule[-2.5truemm]{5cm}{0.1truemm}\\[2mm] {\footnotesize
}
\end{figure}

\section{Introduction and the main results}
\bigskip
\noindent

A {\it general Dirichlet series} is an infinite series of the form
\begin{equation}\label{series}
\sum\limits_{n=1}^{\infty}a_n e^{-\lambda_n s},
\end{equation}
where $a_n$'s are complex numbers, $\{\lambda_n\}$ is a strictly increasing sequence of real numbers whose limit is infinity and $s=\sigma+it$ is a complex variable. The special type $\sum\limits_{n=1}^{\infty}\frac{a_n}{n^s}$, also referred to as {\it ordinary Dirichlet series}, is obtained by substituting $\lambda_n=\log n$. We refer to the monographs \cite{Ap}, \cite{HR}, \cite{Ma}, etc. for general theory and applications of Dirichlet series. Dirichlet series appearing in many problems admit analytic continuation in the complex plane $\textbf{C}$ as meromorphic functions and information on their analytic properties are important in many areas from analysis to number theory to mathematical physics. This paper is concerned with the number of zeros and poles (as well as $a$-values) of these general Dirichlet series. Denote by $n(r, a; f)$ the number of solutions, counted with multiplicity, of the equation $f(s)=a$ in the disc $|s|\le r$ for a function $f$ meromorphic in $\textbf{C}$ and a value $a$. Since we are dealing with general Dirichlet series without assumptions such as functional equation, this {\it counting function} (instead of the one in strips), or its integrated form (see below), is a natural measure for the zeros of $f(s)-a$. This measure, as well-known, is one of the most important quantities in the study of value distribution of meromorphic functions and their applications. The growth of $n(r, a; f)$ for a meromorphic function $f$ can be at the same rate as, but can not exceed, that of the function (see e.g. \cite{Ha}, p.26), which provides a general connection between upper bound of the counting function and growth of the function; there is however no such kind of connection in general between lower bound of the counting function and growth of the function. These bounds have long been studied for Dirichlet series, especially, of course, for the Riemann zeta function $\sum\limits_{n=1}^{\infty}\frac{1}{n^s}$, for which the well-known Riemann-von Mangoldt formula (see e.g. \cite{T2}, p.214) characterizes the counting function of its zeros, and also for other ordinary and general Dirichlet series under various analytic and arithmetic conditions (see e.g. \cite{Le}, \cite{Ra}, \cite{Ra1}, \cite{Se}, \cite{St}, and \cite{T1}, to list a few). It appears in many known cases of Dirichlet series with ``few" poles, besides various L-functions, that $n(r, 0; f)$ grows at least at the order $r$ and there are nontrivial Dirichlet series with growth order equal to $r$. It is thus natural to seek a lower bound of growth order equal to $r$ for the zeros and poles in a rather general setting, and it, if established, would be best possible. In fact, Bombieri and Perelli \cite{BP2} obtained a lower bound, which is of growth order $r^{\gamma}$ for any $\gamma<1$, on the number of zeros and poles for general Dirichlet series under mild analytic assumptions; more precisely, they proved the following

\bigskip
\noindent
{\bf Theorem A.} {\it Let $f(s)=\sum a_ne^{\lambda _n s}, \lambda_n\in \mathbb{R},$ be uniformly convergent in a
half-plane $\Re(s)>b$ and admit an analytic continuation in $\textbf{C}$ as a nonconstant meromorphic function of finite order.
Suppose that $f(s)$ tends to a nonzero finite limit as $\Re(s)\to +\infty$. Then, for any fixed $\gamma<1$, we have
\begin{equation}\label{BP}
\overbar{\lim\limits_{r\to +\infty}}\frac{N(r, 0; f)+N(r, \infty; f)}{r^{\gamma}}>0.
\end{equation}}

\bigskip
\noindent
On the above,
\begin{equation*}
N(r, a; f)=\int_0^r\frac{n(t, a; f)-n(0, a; f)}{t}dt+n(0, a; f)\log r,
\end{equation*}
which is the integrated counting function and plays a similar role as $n(r, a; f)$ for the solutions of the equation $f(s)=a$. (It will be seen from (\ref{un-integrated}) in \S 2 that $N(r, a; f)$ can be replaced by the simpler un-integrated counting function $n(r, a; f)$.)

\medskip
The lower bound given in (\ref{BP}) is not of growth order equal to $r$ due to the assumption $\gamma<1$. Bombieri and Perelli asked whether the conclusion of Theorem A holds with $\gamma=1$, which would be best possible; that is,

\bigskip
\noindent
{\bf Problem A}(\cite{BP2}, p.71). {\it Under the same conditions of Theorem A, whether do we have
\begin{equation}\label{conjecture}
\overbar{\lim\limits_{r\to +\infty}}\frac{N(r, 0; f)+N(r, \infty; f)}{r}>0?
\end{equation}
}

Bombieri and Perelli proved that (\ref{conjecture}) does hold under the following additional conditions \cite{BP2}:

\bigskip
\noindent
{\bf Theorem B}. {\it In addition to the hypotheses of Theorem A, suppose that $f(s)$ is uniformly almost periodic in some half plan $\Re(s)<c$ and $f(s)$ tends to a nonzero finite limit as $\Re (s)\to-\infty.$ Then, (\ref{conjecture}) holds.
}

\bigskip
\noindent
Here, the notion of uniformly almost periodic function is in the sense of Bohr (see \S 2 for the definition).

\bigskip

In this paper, we will solve Problem A by an affirmative answer (and thus the additional conditions in Theorem B are superfluous); that is, we will prove

\bigskip
\noindent
{\bf Theorem 1.1.}  {\it Under the same conditions of Theorem A, we have
\begin{equation*}
\overbar{\lim\limits_{r\to +\infty}}\frac{N(r, 0; f)+N(r, \infty; f)}{r}>0.
\end{equation*}
}

\bigskip
\noindent
{\bf Remark.} We show in this remark that Theorem 1.1 is best possible in the following senses.
\vskip.1in
\noindent
(i) As mentioned above, the growth order $r$ in the conclusion of Theorem 1.1 is the best possible, as seen from the Dirichlet series $f(s)=\sum_{n=0}^{\infty}\frac{a_n}{2^{ns}}$ with suitable $a_n$ (for example, simply take $f(s)=1+\frac{1}{2^{s}}$).
\vskip.1in
\noindent
(ii) The condition `` $f(s)$ tends to a nonzero finite limit as $\Re(s)\to +\infty$" cannot be dropped, as seen from the function $f(s)=e^{-s}$ (resp. $f(s)=e^s$), which tends to zero as $\Re(s)\to +\infty$ (resp. tends to $\infty$ as $\Re(s)\to +\infty$). The function $f$ satisfies all the other conditions of Theorem 1.1 but has no zeros or poles.
\vskip.1in
\noindent
(iii) The condition ``finite order" cannot be dropped in Theorem 1.1. To see this, consider the Dirichlet series
$
f(s)=\sum_{n=0}^{\infty}\frac{1}{n!}e^{-ns},
$
which clearly satisfies all the conditions of Theorem 1.1 except that $f(s)$ is of infinite order, in view of the fact that $f(s)=e^{e^{-s}},$ which has no zeros or poles.
\vskip.1in
\noindent
(iv) The sum $N(r, 0; f)+N(r, \infty; f)$ in the conclusion of Theorem 1.1 cannot be strengthened as $N(r, 0; f)$ or $N(r, \infty; f)$ alone, as seen from the Dirichlet series $f(s)=\frac{1}{1-e^{-s}}=\sum_{n=0}^{\infty}e^{-ns}$ ($\Re (s)>0$), which has no zeros, and the Dirichlet series $f(s)=1-e^{-s}$, which has no poles.

\bigskip
Under the same conditions of Theorem 1.1 or Theorem A, we will actually establish the following result, which improves the original estimate (\ref{conjecture}) in Problem A; specifically we have the following

\bigskip
\noindent
{\bf Theorem 1.2.} {\it Let $f(s)=\sum a_ne^{\lambda _n s}, \lambda_n\in \mathbb{R},$ be uniformly convergent in
a half-plane $\Re(s)>b$ and admit an analytic continuation in $\textbf{C}$ as a nonconstant meromorphic function of finite order.
Suppose that $f(s)$ tends to a nonzero limit as $\Re(s)\to +\infty$. Then either
\begin{equation}\label{conclusion2}
\frac{n(r, 0; f)+n(r, \infty; f)}{r}>A
\end{equation}
for all large $r$, where $A>0$ is a constant, or
\begin{equation}\label{conclusion3}
\int^{\infty}_{r_0}\frac{n(t, 0; f)+n(t, \infty; f)}{t^3}dt=\infty
\end{equation}
for some $r_0>0$.}

\section{Consequences of Theorem 1.2}
\bigskip
\noindent

Theorem 1.2 has various consequences. We first show that Theorem 1.1 is a direct consequence of Theorem 1.2, from which it is also clear (see (\ref{conclusion2}) and (\ref{conclusion4}) below) that {\it under the same conditions of Theorem 1.1 or Theorem 1.2,
\begin{equation}\label{un-integrated}
\overbar{\lim\limits_{r\to +\infty}}\frac{n(r, 0; f)+n(r, \infty; f)}{r}>0.
\end{equation}
}

\medskip
In fact, each of (\ref{conclusion2}) and (\ref{conclusion3}), given in terms of the un-integrated counting function (often more convenient as a counting function), yields a stronger result than the conclusion of Theorem 1.1, i.e., (\ref{conjecture}):
\medskip

If (\ref{conclusion2}) holds, then by the definition of the integrated counting function, we clearly have, by integration, that for some constant $A>0$,
\begin{equation*}
\frac{N(r, 0; f)+N(r, \infty; f)}{r}>A
\end{equation*}
for all large $r$ globally, which is a stronger statement than (\ref{conjecture}).

If (\ref{conclusion3}) holds, then it is easy to see that for any small $\epsilon>0$,
\begin{equation}\label{conclusion4}
\overbar{\lim\limits}_{r\to +\infty}\frac{n(r, 0; f)+n(r, \infty; f)}{r^{2-\epsilon}}=+\infty
\end{equation}
since otherwise there would be an $M>0$ such that
\begin{equation*}
n(r, 0; f)+n(r, \infty; f)<M r^{2-\epsilon}
\end{equation*}
for all large $r$, which clearly implies that
\begin{equation*}
\int^{\infty}\frac{n(t, 0; f)+n(t, \infty; f)}{t^3}dt<+\infty
\end{equation*}
a contradiction to (\ref{conclusion3}). In view of the fact that for $r\ge e$,
\begin{eqnarray*}
& &n(r, 0; f)+n(r, \infty; f)   \nonumber\\
& &\le \int_r^{er}\frac{n(t, 0; f)+n(t, \infty; f)}{t}dt  \nonumber\\
& &\le N(er, 0; f)+N(er, \infty; f),
\end{eqnarray*}
we deduce from (\ref{conclusion4}) the following stronger result than (\ref{conjecture}):
\begin{equation*}
\overbar{\lim\limits}_{r\to +\infty}\frac{N(r, 0; f)+N(r, \infty; f)}{r^{2-\epsilon}}=+\infty.
\end{equation*}

\medskip
In Remark (ii) of Theorem 1.1, we pointed out that the condition `` $f(s)$ tends to a nonzero finite limit as $\Re(s)\to +\infty$" cannot be dropped in Theorem 1.1 and thus in Theorem 1.2. However, this condition can be dropped for general Dirichlet series in both Theorem 1.1 and Theorem 1.2 (Note that the series $f$ in Theorem 1.1 and Theorem 1.2 is an infinite exponential sum, which is more general than general Dirichlet series), except the trivial case that $f(s)=ae^{-\lambda s}$ ($a\not=0$) with only one term, which obviously has no zeros or poles. We particularly include the following corollary, in which (and in the sequel) $f$ is called {\it nontrivial} if $f(s)\not\equiv ae^{-\lambda s}$ with $a, \lambda$ being constant.

\bigskip
\noindent
{\bf Corollary 2.1.}  {\it Let $f$ be a nontrivial general Dirichlet series uniformly convergent in a half-plane $\Re(s)>b$ and admit an analytic continuation in $\textbf{C}$ as a meromorphic function of finite order. Then,
\begin{equation*}
\overbar{\lim\limits_{r\to +\infty}}\frac{n(r, 0; f)+n(r, \infty; f)}{r}>0.
\end{equation*}
In particular, if $f$ has ``few" poles in the sense that $n(r, \infty; f)=o(r)$, then
\begin{equation*}
\overbar{\lim\limits_{r\to +\infty}}\frac{n(r, 0; f)}{r}>0.
\end{equation*}
}

To show Corollary 2.1, suppose that $f(s)=\sum_{n=1}^{\infty}a_n e^{-\lambda_n s}$ is a general Dirichlet series, i.e., $\{\lambda_n\}$ is a sequence of real increasing numbers tending to infinity (see (\ref{series})). We may then assume, without loss of generality, that $a_1\not=0$ (otherwise just use the first nonzero coefficient $a_n$). The Dirichlet series $e^{\lambda_1 s}f(s)=a_1+\sum_{n=2}^{\infty}a_n e^{-(\lambda_n-\lambda_1)s}$ and $f(s)$ have the same zeros and poles. Thus, it suffices to consider  $e^{\lambda_1 s}f(s)$ in Theorem 1.2, which is nonconstant since $f$ is nontrivial, tends to a nonzero finite limit (equal to $a_1$) as $\Re(s)\to +\infty$, and thus satisfies all the conditions of Theorem 1.2. Corollary 2.1 then follows from Theorem 1.2 and (\ref{un-integrated}).

\medskip
The above results readily extend to those for $a$-values of $f$. In particular, for Corollary 2.1, all the conditions are still satisfied for $f-a$, provided that $f-a$ is nontrivial, i.e., $f-a$ is not of the form $ce^{-\lambda s}$ with $c, \lambda$ being constant (cf. Corollary 2.1). Thus, we obtain from Corollary 2.1 the following quantitative Picard type theorem:

\bigskip
\noindent
{\bf Corollary 2.2.}  {\it Let $f$ be a general Dirichlet series uniformly convergent in a half-plane $\Re(s)>b$ and admit an analytic continuation in $\textbf{C}$ as a meromorphic function of finite order. Then
\begin{equation}\label{a-values}
\overbar{\lim\limits_{r\to +\infty}}\frac{n(r, a; f)+n(r, \infty; f)}{r}>0
\end{equation}
for any complex number $a$, provided that $f-a$ is nontrivial.
}

\medskip
It is worth noting that in Corollary 2.2, if $f$ is not of the form $a_0+b_0e^{-\lambda s}$ for some constants $\lambda, a_0, b_0$, then (\ref{a-values}) holds for every complex number $a$ without any exception (The situation when $f$ is of the form $a_0+b_0e^{-\lambda s}$ is trivial); in the case that $f$ is entire or has ``few" poles in the sense stated above, i.e., $n(r, \infty; f)=o(r)$, the conclusion (\ref{a-values}), giving a lower bound on how often $f$ assumes $a$, implies particularly that $f$ takes every complex number $a$ infinitely often without any exception. This should be compared to the famous Picard theorem: A nonconstant entire (resp. meromorphic) function takes every complex number infinitely often with possibly one exception (resp. two exceptions), and can be thought of as a quantitative Picard type Theorem for general Dirichlet series. We refer to the monograph \cite{Ma}, Chapter III, for related but different results for Picard type theorems for Dirichlet series.

\bigskip
In many problems of significance in number theory, Dirichlet series admit an analytic continuation in $\textbf{C}$ as meromorphic functions of order $\le 1$. Theorem 1.2 may further be applied to yield a global estimate, which is substantially stronger than the estimate from the upper limit result, on the counting functions of zeros, poles, and, more generally, $a$-values for these Dirichlet series $f$. In fact, when $f$ is of order $\le 1$, it follows from Jensen's formula (see e.g. \cite{Ha}, p.25-26) that for any complex value $a$ (including $a=\infty$), $n(r, a; f)=O(r^{1+\epsilon})$ for any $0<\epsilon<1$, which implies that
\begin{equation*}
\int^{\infty}\frac{n(t, 0; f)+n(t, \infty; f)}{t^3}dt<+\infty.
\end{equation*}
This last inequality holds actually for all $f$ of order less than $2$, since then $n(r, a; f)=O(r^{2-\delta})$ for some small $\delta>0$. Thus for these $f$, (\ref{conclusion3}) does not hold, and consequently, the global estimate (\ref{conclusion2}) in Theorem 1.2 must hold. That is, we have the following (cf. Corollary 2.1)

\bigskip
\noindent
{\bf Corollary 2.3.}  {\it Let $f$ be a nontrivial general Dirichlet series uniformly convergent in a half-plane $\Re(s)>b$ and admit an analytic continuation in $\textbf{C}$ as a meromorphic function of order $<2$. Then
\begin{equation}\label{global}
n(r, 0; f)+n(r, \infty; f)>Ar
\end{equation}
for all large $r$, where $A>0$ is a constant.
}

\bigskip
As in Theorem 1.2, the nontrivial general Dirichlet series $f$ in Corollary 2.3 can be replaced by a nonconstant infinite exponential sum $f(s)=\sum a_ne^{\lambda _n s}, \lambda_n\in \mathbb{R}$ with $f(s)$ tending to a nonzero limit as $\Re(s)\to +\infty$, and the result also extends to $a$-values of $f$ as given in Corollary 2.2. We omit their statements here.

\bigskip
We next present another consequence of Theorem 1.2 (Corollary 2.3) as an application to L-functions, which are extensively studied in number theory, for lower bounds on the cardinality of the symmetric difference of their zeros. Let $F(s)$ and $G(s)$ be two general Dirichlet series, and $m_F(\rho)$ (resp. $m_G(\rho)$) the multiplicity of $\rho$ as zero of $F$ (resp. of $G$). The symmetric difference $D_{F, G}$ is defined as
\begin{equation}\label{difference}
D_{F, G}=\sum_{|\rho|\le T}|m_F(\rho)-m_G(\rho)|,
\end{equation}
where $\rho=\sigma+it$ runs over the zeros of $F(s)G(s)$. As given in \cite{BP2}, the following holds on the cardinality $D_{F, G}$ of two distinct L-functions $F, G$ satisfying the same functional equation:
\begin{equation}\label{symmetry}
D_{F, G}(T)=\Omega(T)
\end{equation}
(counting only the nontrivial zeros of $F(s)G(s)$ with $|t|\le T$ in (\ref{difference})), where the notation $f(x)=\Omega(g(x)$ denotes the negation of $f(x)=o(g(x))$ (see e.g. \cite{T2}, p.184), or $\overbar{\lim\limits}_{T\to +\infty}\frac{D_{F, G}(T)}{T}>0$, under only the function-theoretic properties of $F$ and $G$, disregarding their arithmetical aspects. Here and in the sequel, by abuse of language, an {\it L-function} means a Dirichlet series $F(s)=\sum\limits_{n=1}^{\infty}\frac{a_n}{n^s}$ with $a_1=1$ which is absolutely convergent for $\sigma>1$, admits an analytic continuation in $\textbf{C}$ as a meromorphic function of finite order with at most one pole at $s=1$, and satisfies a standard functional equation (cf. \cite{Se}), without requiring the Ramanujan hypothesis and Euler product. Note that the same bound (\ref{symmetry}) was obtained in \cite{MM} for any two distinct L-functions $F$ and $G$ in the Selberg class, which consists of Dirichlet series $\sum_{n=1}^{\infty}\frac{a_n}{n^s}$ satisfying the above mentioned conditions and also Ramanujan hypothesis and Euler product(see \cite{Se} for the details). Now, by the fact that a nonconstant L-function, as defined above, is actually of order equal to $1$ (see e.g. \cite{Se}, p.48), which implies that for two L-functions $F$ and $G$, the Dirichlet series $f=\frac{F}{G}$ (not an L-function in general) is of order at most $1$ and thus from Corollary 2.3, the inequality (\ref{global}) holds for $f=\frac{F}{G}$, which clearly improves the above mentioned result (\ref{symmetry}) of \cite{BP2} to the following {\it global estimate for all large $T$:
\begin{equation}\label{global-1}
D_{F, G}(T)>AT,
\end{equation}
where $A>0$ is a constant.}

Note from the above that the estimate (\ref{global-1}) holds also for general Dirichlet series under much weaker conditions than those for L-functions. It also connects to another topic- uniqueness problems for L-functions, or more generally, Dirichlet series, about when $F$ and $G$ are identically equal in terms of the so-called shared values, which has recently been studied in various settings, see \cite{CY}, \cite{GL}, \cite{GHK}, \cite{HL}, \cite{Ki}, \cite{KL}, \cite{Li}, \cite{St}, etc. From Corollary 2.3 and (\ref{global-1}), we have the following

\bigskip
\noindent
{\bf Corollary 2.4.}  {\it Let $F$ and $G$ be two L-functions, or, more generally, let $F, G$ be two general Dirichlet series so that $F/G$ is uniformly convergent in a half-plane $\Re(s)>b$ and admit an analytic continuation in $\textbf{C}$ as a meromorphic function of order $<2$ with the limit equal to $1$ as $\Re(s)\to +\infty$. Then $F\not\equiv G$ if and only if
\begin{equation}\label{global-2}
D_{F, G}(T)>AT \quad \hbox{or}\quad  n(T, 0; F/G)+n(T, \infty; F/G)>AT
\end{equation}
for all large $T$, where $A>0$ is a constant. Equivalently, $F\equiv G$ if and only if (\ref{global-2}) does not hold.}

\bigskip
The assumption in Corollary 2.4 that $F/G$ tends to $1$ as $\Re(s)\to +\infty$ together with $F\not\equiv G$ guarantees that $F/G$ is nontrivial; otherwise, $F/G=ae^{-\lambda s}$ for some $a, \lambda$, from which it follows that $\lambda=0$ and then $a=1$ by letting $\Re(s)\to +\infty$ and thus $F\equiv G$, a contradiction. It is interesting to see that this limit condition is tied to the exceptional case in a more familiar form of uniqueness in terms of shared values (cf. Corollary 2.5 below).

If two L-functions $F$ and $G$ share the value $0$, i.e., $F$ and $G$ have the same zeros with counting multiplicities, then it is known that $F\equiv G$ (see \cite{St}, p.152), which can now follow from Corollary 2.4 immediately since $D_{F, G}(T)=O(1)$ and thus (\ref{global-2}) does not hold and consequently $F\equiv G$. If two L-functions $F$ and $G$ share a complex number $a$, i.e., $F-a$ and $G-a$ have the same zeros with counting multiplicities, we may apply Corollary 2.4 to $F-a$ and $G-a$. Note that $\frac{F-a}{G-a}$ has the limit $1$ as $\Re(s)\to +\infty$, except when $a=1$, for which the limit might not be $1$. We deduce from Corollary 2.4 the uniqueness that $F-a\equiv G-a$, i.e., $F\equiv G$, except the case $a=1$, for which the result need not hold (compare the result in \cite{St}, p.152; see also \cite{HL}). In fact, the quantitative condition (\ref{global-2}) can produce a quantitative  condition for uniqueness, see Corollary 2.5 below, in which we say that $F-a$ and $G-a$ {\it have ``enough" common zeros (resp. poles)} if $F-a$ and $G-a$ have the same zeros (resp. poles) counted with multiplicity except an exceptional set $E$ of their zeros (resp. poles) satisfying that $n(r, E)=o\{r\}$ as $r\to\infty$, where $n(r, E)$ is the counting function of $E$, i.e., the number of points in $E\cap \{|s|<r\}$ counted with multiplicity. From Corollary 2.4, we immediately have the following uniqueness theorem.

\bigskip
\noindent
{\bf Corollary 2.5.}  {\it Let $F$ and $G$ be two general Dirichlet series and $a\not=1$ a complex number so that $\frac{F-a}{G-a}$ is uniformly convergent in a half-plane $\Re(s)>b$ and admit an analytic continuation in $\textbf{C}$ as a meromorphic function of order $<2$ with the limit equal to $1$ as $\Re(s)\to +\infty$. Then $F\equiv G$ if and only if $F-a$ and $G-a$ have ``enough" common zeros and poles.

In particular, if $F$ and $G$ are two L-functions, then $F\equiv G$ if and only if $F-a$ and $G-a$ have ``enough" common zeros.

The result need not hold when $a=1$.}

\bigskip
To see Corollary 2.5 need not hold when $a=1$, consider the functions $F=1+\frac{2}{4^s}$ and $G= 1+\frac{3}{9^s}$. Then $F, G$ satisfy trivially all the sufficient conditions of Corollary 2.5 with $a=1$. Actually, $F$ and $G$ satisfy the functional equations $2^sL(s)=2^{1-s}\overline{L(1-{\bar
s})}$ and
$3^sL(s)=3^{1-s}\overline{L(1-{\bar s})}$, respectively (cf. \cite{HL}) and thus it is easy to see that $F$ and $G$ are also L-functions as defined above. But, $F\not\equiv G$.

\medskip
We note that a connection (an equivalence) between the Riemann hypothesis and a uniqueness problem for the Riemann zeta function $\zeta(s)$ and its analogue for L-functions has been given in \cite{HL2}.

\section{Proof of Theorem 1.2}
\bigskip
\noindent

Throughout the proof, we will use $\epsilon$ to denote a positive constant which can be made arbitrarily small and $C$ a positive constant, the actual values of which may vary from one occurrence to the next.

Recall the classic Poisson-Jensen formula (see e.g. \cite{La}, p.345 and \cite{Ha}, p.1): Let $f$ be a meromorphic function in $|s|\le R$. Then for each $s$ in $|s|<R$ with $f(s)\not=0, \infty$, we have
\begin{eqnarray}\label{poisson}
& &\log |f(s)|=\frac{1}{2\pi}\int_0^{2\pi}\Re\{\frac{Re^{i\phi}+s}{Re^{i\phi}-s}\}\log|f(Re^{i\phi})|d\phi+\\
& &+\sum_{k=1}^{M}\log|\frac{R(s-a_k)}{R^2-\bar{a_k}s}|-\sum_{k=1}^{N}\log|\frac{R(s-b_k)}{R^2-\bar{b_k}s}|, \nonumber
\end{eqnarray}
where $a_k (1\le k\le M)$ and $b_k (1\le k\le N)$ are the zeros and poles of $f$ in $|s|<R,$ respectively. In particular, if $s=0$, we have the Jensen formula:
\begin{eqnarray}\label{jensen}
\log |f(0)|=\frac{1}{2\pi}\int_0^{2\pi}\log|f(Re^{i\phi})|d\phi+\sum_{k=1}^{M}\log|\frac{a_k}{R}|-\sum_{k=1}^{N}\log|\frac{b_k}{R}|.
\end{eqnarray}
If $f(s)=0$ or $\infty$, $\log |f(0)|$ is then replaced by $\log a_m+m\log R,$ where $m=n(0, 0; f)-n(0, \infty; f)$ and $a_m\not=0$ is the coefficient in the Laurant or Taylor expansion $f=a_ms^m+\cdots$ of $f$ at $s=0$.

The underlying ideas of utilizing the Poisson-Jensen formula, which connects the modulus (and thus the growth) of a meromorphic function with its zeros and poles, can be traced to the relevant techniques in Nevanlinna theory (see \cite{Ha}, p.36, \cite{DDL}, \cite{CF}, etc.), while the present paper does not need familiarity with the theory. The Poisson-Jensen formula will be combined with several other tools in our proof.

We will make use of the following result (see \cite{Ha}, p.27 and \cite{GO}, p.56), which can be proved directly by estimating the modulus of the infinite product below: Let $\{w_n\}$ be a sequence of nonzero complex numbers such that $\sum |w_n|^{-2}$ converges. Then
the infinite product $\Pi(s):=\Pi_{n=1}^{\infty}(1-\frac{s}{w_n})e^{\frac{s}{w_n}}$ is an entire function satisfying the following estimate
\begin{eqnarray}\label{estimate1}
& &\log|\Pi(s)| \\
& &\leq 4(2+\log 2)\{|s|\int_0^{|s|}\frac{n(t, 0, \Pi)}{t^2}dt+|s|^2\int_{|s|}^{+\infty}\frac{n(t, 0, \Pi)}{t^3}dt\}. \nonumber
\end{eqnarray}

We will also employ the following result due to Cartan (see e.g. \cite{BG}, p.360 and \cite{Lev}, p.19): Given any number $h>0$ and complex numbers $a_1, a_2,\cdots, a_n$, there is a collection of disks in the complex plane with the sum of the radii equal to $2h$ such that for each point $s$ lying outside these disks one has the inequality
\begin{eqnarray}\label{cartan}
\Pi_{k=1}^n|s-a_k|>(\frac{h}{e})^n.
\end{eqnarray}

Recall also that an analytic function $f(s)$ is uniformly almost periodic (in the sense of Bohr) in a strip $b < \Re(s)< c $ ($b$ and $c$ may be $\pm \infty$) if for every $\epsilon>0$, the set of real numbers $\omega$ such that
\begin{equation*}
|f (s + i\omega ) - f (s)| < \epsilon,   \quad b <\Re(s)< c
\end{equation*}
is relatively dense, i.e., if for every $\epsilon>0$ there is an $l>0$ such that every interval of length $l$ contains such a number $\omega$ (see e.g. \cite{Be}, \cite{Bo}). If $f$ is a nonconstant uniformly almost periodic function in a strip $b < \Re(s)< c $ and $f(s)=a$ has a solution in the strip, then we have that
\begin{equation}\label{counting}
N(r, a; f)\ge cr
\end{equation}
for some constant $c>0$ and all large $r$. This follows from Rouch\'e's theorem in a standard argument (see e.g. \cite{BP2}, \cite{DH}). For completeness, we include its proof. In fact, suppose that $s_0$ is a solution of $f(s)=a$ in the strip and let $\epsilon$ be the minimum of $|f(s)-a|$ on the boundary of a disk containing $s_0$ and lying inside the strip with $f(s)\not=a$ on the boundary, then by uniformly almost periodic property, there is a number $l>0$ such that every interval of length $l$ contains a $\tau$ such that $|f(s+i\tau)-f(s)|<\epsilon$ on the boundary. By Rouch\'e's theorem, $f(s)-a$ and $f(s+i\tau)-a$ have the same number of zeros inside the disc, which clearly implies the inequality (\ref{counting}).

\bigskip
Let us now turn to the proof of Theorem 1.2.

\bigskip
\noindent
{\bf Proof of Theorem 1.2}. If the conclusion (\ref{conclusion3}) holds, we have nothing to prove. Thus, it suffices to prove (\ref{conclusion2}) when
\begin{equation}\label{con3}
\int^{\infty}_{r_0}\frac{n(t, 0; f)+n(t, \infty; f)}{t^3}dt<\infty
\end{equation}
for some $r_0>0$.

Suppose, to the contrary, that (\ref{conclusion2}) does not hold. Then there exists an increasing sequence $\{r_n\}$ of positive numbers with $r_n\to +\infty$ such that for all $n$,
\begin{equation}\label{2.1}
\frac{n(r_n, 0; f)+n(r_n, \infty; f)}{r_n}< \frac{1}{n}.
\end{equation}
We will arrive at a contradiction eventually.

By (\ref{con3}), for any given $\epsilon>0$ we have that for large $r$,
\begin{eqnarray*}
& &\epsilon>\int_r^{2r}\frac{n(t, 0; f)+n(t, \infty; f)}{t^3}dt\\
& &\geq \frac{n(r, 0; f)+n(r, \infty; f)}{(2r)^3}r\\
& &=\frac{1}{8}\frac{n(r, 0; f)+n(r, \infty; f)}{r^2},
\end{eqnarray*}
which implies that
\begin{equation}\label{limit}
\frac{n(r, 0; f)+n(r, \infty; f)}{r^2}\to 0
\end{equation}
as $r\to +\infty$.

Let  $a_k$ ($k=1, 2,\cdots$) be the nonzero zeros of $f$ arranged with $|a_k|\le |a_{k+1}|$ and $b_k$ ($k=1, 2,\cdots$) the nonzero poles of $f$ arranged with $|b_k|\le |b_{k+1}|$. Then by the Stieltjes integral and (\ref{limit}),
\begin{eqnarray*}
& &\sum_{k=1}^{\infty}|a_k|^{-2}=\int_0^{\infty}\frac{1}{t^2}{d(n(t, 0; f)-n(0, 0; f)}\\
& &=\lim\limits_{r\to \infty}(\frac{n(r, 0; f)-n(0, 0; f)}{r^2}+2\int_0^r\frac{n(t, 0; f)-n(0, 0; f)}{t^3}dt)\\
& &=2\int_0^{\infty}\frac{n(t, 0; f)-n(0, 0; f)}{t^3}dt
\end{eqnarray*}
converges by (\ref{con3}). One can then check (or see (\ref{estimate1})) that the infinite product
\begin{equation*}
f_1=\Pi_{k=1}^{\infty}(1-\frac{s}{a_k})e^{\frac{s}{a_k}}
\end{equation*}
is an entire function having $a_k$'s as its zeros. Similarly, the infinite product
\begin{equation*}
f_2=\Pi_{k=1}^{\infty}(1-\frac{s}{b_k})e^{\frac{s}{b_k}}
\end{equation*}
is an entire function having $b_k$'s as its zeros. Therefore, we have the following factorization of $f$:
\begin{equation*}
f(s)=s^m\frac{f_1(s)}{f_2(s)}e^{Q(s)},
\end{equation*}
where $m$ is an integer, which is the order of zero or pole of $f$ at $s=0$, and $Q$ is a polynomial of degree $\deg(Q)\le \rho(f)$, where $\rho(f)$ denotes the order of $f$.

Fix a small positive number $\tau$ (to be specified later). Define the operator $\Lambda$:
\begin{equation*}
\Lambda f(s)=\frac{f(s+\tau)}{f(s)}, \quad \Lambda^k=\Lambda(\Lambda^{k-1}).
\end{equation*}
It is clear that  $\Lambda e^{Q(s)}=e^{Q(z+\tau)-Q(z)}$ and $Q(z+\tau)-Q(z)$ is a polynomial of degree at most $\deg(Q)-1$. Applying the operator  $\Lambda$ again, we see that $\Lambda^2 e^{Q(s)}=e^{Q_1(s)}$, where $Q_1(s)$ is a polynomial of degree at most $\deg(Q)-2$. Applying $\Lambda$ repeatedly for $d_0$ times, we obtain that $\Lambda^{d_0}e^{Q(s)}=c_0$, a nonzero constant, where $d_0=\deg(Q)$. Set $d=d_0+1\geq 1$. We thus obtain, in view of the above factorization of $f$,  that
\begin{equation}\label{functionF}
F(s):=\Lambda^{d}f(s)=\frac{\Lambda^{d}f_1(s)}{\Lambda^{d}f_2(s)}R(s),
\end{equation}
where $R(s)=\Lambda^{d}s^m$ is a nonzero rational function. (It is worth mentioning that $d\geq 1$ so we apply the operator $\Lambda$ on $f$ at least once even when the polynomial $Q$ is a constant.)

By the given condition, $f(s)$ is nonconstant and tends to a nonzero finite limit, denoted by $A_0$, as $\Re(s)\to +\infty$. Thus, $\Lambda f(s)$ tends to $1$ as $\Re(s)\to +\infty$. This implies that $\Lambda f(s)$ cannot be a constant function, since otherwise it must be identically equal to $1$; that is, $f(s+\tau)\equiv f(s)$, a periodic function. But, $f$ is bounded when $\Re(s)\to +\infty$. Hence, $f$ is bounded in the complex plane (For any $\sigma_0+it_0$, $|f(\sigma+it_0)|\le |A_0|+1$ for large $\sigma$ and thus for all $\sigma$, particularly $\sigma_0$, by the periodicity) and consequently must be a constant by Liouville's theorem, a contradiction. This same argument (cf. \cite{BP2}) shows that $\Lambda^2 f(s), \cdots, \Lambda^d f(s)$ are all nonconstant. In particular, $F(s)=\Lambda^{d}f(s)$ is nonconstant. Note that an exponential series of the form $\sum_{n=1}^{\infty}a_n e^{\lambda_n s}$ uniformly convergent in a strip $b<\Re(s)<c$ ($b$ and $c$ may be $\pm \infty$) must be uniformly almost periodic in the strip (see e.g. \cite{Be}, Theorem 6, Ch. III). Thus, the Dirichlet series $f$ is uniformly almost periodic in the half-plane $\Re(s)>b$. It is then easy to verify directly that $\Lambda f$, as the quotient of $f(s+\tau)$ and $f(s)$, which tend to a nonzero finite limit $A_0$ as $\Re(s)\to +\infty$, is also uniformly almost periodic in a right half-plane. In fact, we have that $\frac{|A_0|}{2}<|f(s)|<2|A_0|$ and $\frac{|A_0|}{2}<|f(s+\tau)|<2|A_0|$ for large $\Re(s)$. For any $\epsilon>0$, there is an $l>0$ such that every interval of length $l$ contains a $\omega$ such that $|f(s+i\omega)-f(s)|<\epsilon$ for $\Re(s)>b$. Thus,
\begin{eqnarray*}
& &|\Lambda f(s+i\omega)-\Lambda f(s)|\\
& &=|\frac{f(s+\tau+i\omega)}{f(s+i\omega)}-\frac{f(s+\tau)}{f(s)}|\\
& &=|\frac{\bigl(f(s+\tau+i\omega)-f(s+\tau)\bigr)f(s)+f(s+\tau)(f(s)-f(s+i\omega))}{f(s+i\omega)f(s)}|\\
& &\le \frac{4|A_0|\epsilon}{(\frac {|A_0|}{2})^2}=\frac{16\epsilon}{|A_0|}
\end{eqnarray*}
for large $\Re(s)$. This shows that $\Lambda f(s)$ is uniformly almost periodic in a right half-plane. The same argument shows that $\Lambda^2 f(s), \cdots, \Lambda^d f(s)$ are all uniformly almost periodic in a right half-plane; in particular, $F(s)=\Lambda^d f(s)$ is uniformly almost periodic in a right half-plane. Taking a point $w_0$ in this right half-plane and denote $a=F(w_0)$. Then by virtue of (\ref{counting}), there exists a constant $c>0$ such that for large $r$,
\begin{equation}\label{integrated counting}
N(r, a, F)\ge cr.
\end{equation}

We will establish a tight estimate on $N(r, a, F)$ in terms of the zeros of $f_1$ and $f_2$ and thus the zeros and poles of $f$. Fix any $r>e$. By the Jensen formula (\ref{jensen}), we have that
\begin{eqnarray*}
& &\log |F(0)-a|\\
& &=\frac{1}{2\pi}\int_0^{2\pi}\log|F(re^{i\theta}-a|d\theta+\sum_{|\alpha_k|<r}\log\frac{|\alpha_k|}{r}-\sum_{|\beta_k|<r}\log\frac{|\beta_k|}{r},
\end{eqnarray*}
where $\alpha_k$ are the zeros of $F-a$ and $\beta_k$ are the poles of $F-a$ (and thus $F$) in $|s|<r$. Here, if $F(0)=a$ or $\infty$, the term $\log |F(0)-a|$ needs to be replaced by $c_0+c_1\log r$ for some constants $c_0$ and $c_1$. Noting that
\begin{eqnarray}\label{sum}
& &\sum_{|\alpha_k|<r}\log\frac{r}{|\alpha_k|}=\int_0^r\log \frac{r}{t}d(n(t, a; F)-n(0, a; F)) \nonumber\\
& &=\int_0^r\frac{n(t, a; F)-n(0, a; F)}{t}dt \nonumber\\
& &=N(r, a; F)-n(0, a; F)\log r
\end{eqnarray}
and similarly,
\begin{equation*}
\sum_{|\beta_k|<r}\log\frac{r}{|\beta_k|}=N(r, \infty; F)-n(0, \infty; F)\log r,
\end{equation*}
we deduce that
\begin{eqnarray*}
& &N(r, a; F)\le\frac{1}{2\pi}\int_0^{2\pi}\log|F(re^{i\theta}-a|d\theta+N(r, \infty; F)+C\log r\\
& &=\frac{1}{2\pi}\int_{|F(re^{i\theta})|\leq |a|}\log|F(re^{i\theta}-a|d\theta+\frac{1}{2\pi}\int_{|F(re^{i\theta})|>|a|}\log|F(re^{i\theta}-a|d\theta\\
& &+N(r, \infty; F)+C\log r\\
& &\le \frac{1}{2\pi}\int_{|F(re^{i\theta})|\leq |a|}\log(2|a|)d\theta+\frac{1}{2\pi}\int_{|F(re^{i\theta})|>|a|}\log(2|F(re^{i\theta}|)d\theta\\
& &+N(r, \infty; F)+C\log r\\
& &\le\frac{1}{2\pi}\int_{E}\log|F(re^{i\theta}|d\theta+N(r, \infty; F)+C\log r\\
& &\le\frac{1}{2\pi}\int_{E}\log|\frac{\Lambda^{d}f_1(re^{i\theta})}{\Lambda^{d}f_2(re^{i\theta})}|d\theta+N(r, \infty; F)+C\log r,
\end{eqnarray*}
where $E$ denotes the set of $\theta$ such that $|F(re^{\theta})|>|a|,$ in view of (\ref{functionF}) and the fact that
\begin{eqnarray*}
& &\int_E\log |R(re^{\theta})|d\theta=O(\log r)
\end{eqnarray*}
since $R$ is a rational function and thus $\log |R(re^{\theta})|=O(\log r)$. By the definition of the operator $\Lambda$, it is easy to check that $\frac{\Lambda^{d}f_1(s)}{\Lambda^{d}f_2(s)}$ is a product of finitely
many quotients of the form
\begin{equation*}
\frac{f_1(s+\eta+\tau)}{f_1(s+\eta)}, \frac{f_1(s+\eta)}{f_1(s+\eta+\tau)}, \frac{f_2(s+\eta+\tau)}{f_2(s+\eta)}, \frac{f_2(s+\eta)}{f_2(s+\eta+\tau)},
\end{equation*}
where $\eta=k\tau$ and $k$ is an integer with $0\le k\le d-1$, i.e., there are finite sets $I_j (1\le j\le 4),$ whose elements are of the form $k\tau$ with $0\le k\le d-1$, such that
\begin{eqnarray*}
& &\frac{\Lambda^{d}f_1(s)}{\Lambda^{d}f_2(s)}\\
& &=\Pi_{\eta\in I_1}\frac{f_1(s+\eta+\tau)}{f_1(s+\eta)}\Pi_{\eta\in I_2}\frac{f_1(s+\eta)}{f_1(s+\eta+\tau)}\Pi_{\eta\in I_3}\frac{f_2(s+\eta+\tau)}
{f_2(s+\eta)}\Pi_{\eta\in I_4}\frac{f_2(s+\eta)}{f_2(s+\eta+\tau)}.
\end{eqnarray*}
Thus, we obtain that with $s=re^{i\theta}$,
\begin{eqnarray}\label{2.3}
& &N(r, a, F)\leq \sum_{\eta\in I_1}\frac{1}{2\pi}\int_{E}\log|\frac{f_1(s+\eta+\tau)}{f_1(s+\eta)}|d\theta \nonumber\\
& &+\sum_{\eta\in I_2}\frac{1}{2\pi}\int_{E}\log|\frac{f_1(s+\eta)}{f_1(s+\eta+\tau)}|d\theta+\sum_{\eta\in I_3}\frac{1}{2\pi}\int_{E}\log|\frac{f_2(s+\eta+\tau)}
{f_2(s+\eta)}|d\theta \nonumber\\
& &+\sum_{\eta\in I_4}\frac{1}{2\pi}\int_{E}\log|\frac{f_2(s+\eta)}{f_2(s+\eta+\tau)}|d\theta+N(r, \infty; F)+C\log r.
\end{eqnarray}

We now estimate each term in (\ref{2.3}). To estimate $N(r, \infty; F)$, we note that the nonzero poles of $F$ come from the zeros of
$f_1(s+k\tau)$ or $f_2(s+k\tau),$ where $k$ is an integer with $0\le k\le d$. When $|s|\le r$, all these zeros lie within $|s|\le r+d\tau.$ Thus,
\begin{eqnarray*}
& &n(r, 0; f_1(s+k\tau))\leq n(r+d\tau, 0; f_1(s)), \nonumber \\
& &n(r, 0; f_2(s+k\tau))\leq n(r+d\tau, 0; f_2(s)),
\end{eqnarray*}
which clearly implies that
\begin{equation*}
n(r, \infty, F)\le C\bigl(n(r+d\tau, 0; f_1)+n(r+d\tau, 0; f_2)+1\bigr)
\end{equation*}
for all $r>0$, where $C>0$ is a constant. We may take a small $\tau$ in the beginning so that $f_1, f_2$ do not vanish in $|s|\le 2d\tau$ (such a $\tau$ clearly exists
since $f_1(0)f_2(0)\not=0$ by the definitions of $f_1$ and $f_2$). Then, for a fixed small $\delta>0$,
\begin{eqnarray}\label{2.4}
& &N(r, \infty, F)=\int_0^r\frac{n(t, \infty; F)-n(0, \infty; F)}{t}dt+n(0, \infty; F)\log r \nonumber\\
& &=\int_{\delta}^r\frac{n(t, \infty; F)-n(0, \infty; F)}{t}dt+n(0, \infty; F)\log r \nonumber\\
& &\le \int_{\delta}^r\frac{C(n(t+d\tau, 0; f_1)+n(t+d\tau, 0; f_2)+1)}{t}dt+n(0, \infty; F)\log r \nonumber\\
& &\le C(\int_0^r\frac{n(t+d\tau, 0; f_1)+n(t+d\tau, 0; f_2)}{t}dt+\log r) \nonumber\\
& &=C(\int_{d\tau}^{r+d\tau}\frac{n(t, 0; f_1)+n(t, 0; f_2)}{t-d\tau}dt+\log r)\nonumber\\
& &=C\bigl(\int_{2d\tau}^{r+d\tau}\frac{n(t, 0; f_1)+n(t, 0; f_2)}{t}(1+\frac{d\tau}{t-d\tau})dt+\log r\bigr) \nonumber\\
& &\le 2C\bigl(N(r+d\tau, 0; f_1)+N(r+d\tau, 0; f_2)+\log r\bigr) \nonumber\\
& &\le 2C\bigl(N(r+d\tau, 0; f)+N(r+d\tau, \infty; f)+\log r\bigr).
\end{eqnarray}

Next we estimate $\frac{1}{2\pi}\int_{E}\log|\frac{f_1(s+\eta+\tau)}{f_1(s+\eta)}|d\theta$ in (\ref{2.3}); the same estimate for it given below applies to each of the integrals in(\ref{2.3}). To this end, we use the Poisson-Jensen formula (\ref{poisson}) for the entire function $f_1$ in $|s|<R$, where $R$ is a large number to be determined later, but large enough at this moment so that the considered points $s=re^{i\theta}$ and $s+k\tau$ ($1\le k\le d$) all lie
within $|s|<R,$ i.e., $R>r+d\tau$. For all these $s=re^{i\theta}$, except those with $f_1(s+\eta)=0$ or $f_1(s+\eta+\tau)=0$, which constitute a finite set and do not affect the involved integrals below,
we have that
\begin{eqnarray*}
& &\log|f_1(s+\eta)|\\
& &=\frac{1}{2\pi}\int_0^{2\pi}\Re\{\frac{Re^{i\phi}+s+\eta}{Re^{i\phi}-(s+\eta)}\}\log|f_1(Re^{i\phi})|d\phi+\\
& &+\sum_{k=1}^{M}\log|\frac{R(s+\eta-a_k)}{R^2-\bar{a_k}(s+\eta)}|
\end{eqnarray*}
and that
\begin{eqnarray*}
& &\log|f_1(s+\eta+\tau)|\\
& &=\frac{1}{2\pi}\int_0^{2\pi}\Re\{\frac{Re^{i\phi}+s+\eta+\tau}{Re^{i\phi}-(s+\eta+\tau)}\}\log|f_1(Re^{i\phi})|d\phi+\\
& &+\sum_{k=1}^{M}\log|\frac{R(s+\eta+\tau-a_k)}{R^2-\bar{a_k}(s+\eta+\tau)}|,
\end{eqnarray*}
where $a_k$ ($k=1, 2\cdots, M$) are the zeros of $f_1$ in $|s|<R.$ Observing that
\begin{equation}\label{2.5}
\frac{R^2-\bar{a_k}s}{R(s-a_k)}\ge 1
\end{equation}
for any $|s|\leq R$, we deduce that
\begin{eqnarray}\label{2.6}
& &\log|\frac{f_1(s+\eta+\tau)}{f_1(s+\eta)}| \nonumber \\
& &= \frac{1}{2\pi}\int_0^{2\pi}\Re\{\frac{Re^{i\phi}+s+\eta+\tau}{Re^{i\phi}-(s+\eta+\tau)}-\frac{Re^{i\phi}+s+\eta}{Re^{i\phi}-
(s+\eta)}\}\log|f_1(Re^{i\phi})|d\phi+ \nonumber\\
& &+\sum_{k=1}^{M}\log|\frac{R(s+\eta+\tau-a_k)}{R^2-\bar{a_k}(s+\eta+\tau)}|-\sum_{k=1}^{M}\log|\frac{R(s+\eta-a_k)}{R^2-\bar{a_k}(s+\eta)}| \nonumber\\
& &\le \frac{1}{2\pi}\int_0^{2\pi}\frac{2\tau R}{(R-r-d\tau)^2}|\log|f_1(Re^{i\phi})||d\phi \nonumber\\
& &+\sum_{k=1}^{M}\log|\frac{R^2-\bar{a_k}(s+\eta)}{R(s+\eta-a_k)}|.
\end{eqnarray}
It is easy to check that for any $x>0$, $|\log x|=2\log^+x-\log x$, where $\log^+x=\max\{\log x, 0\}$. Thus,
\begin{eqnarray} \label{2.7}
& &\frac{1}{2\pi}\int_0^{2\pi}|\log|f_1(Re^{i\phi})||d\phi  \nonumber\\
& &=\frac{1}{\pi}\int_0^{2\pi}\log^+|f_1(Re^{i\phi})|d\phi-\frac{1}{2\pi}\int_0^{2\pi}\log |f_1(Re^{i\phi})|d\phi \nonumber\\
& &\le \frac{1}{\pi}\int_0^{2\pi}\log^+|f_1(Re^{i\phi})|-\log|f_1(0)|
\end{eqnarray}
by virtue of the Jensen formula (\ref{jensen}) for the last inequality. We further estimate the last integral in (\ref{2.7}) using the tight estimate provided by (\ref{estimate1}). Applying this estimate to $f_1$ and in view of (\ref{con3}) and (\ref{limit}), we obtain that for any given $\epsilon>0$,
\begin{eqnarray*}
& &\log |f_1(s)|\leq 4(2+\log 2)\{|s|\int_0^{|s|}\frac{n(t, 0, f_1)}{t^2}dt+|s|^2\int_{|s|}^{+\infty}\frac{n(t, 0, f_1)}{t^3}dt\}\\
& &\le \epsilon|s|^2
\end{eqnarray*}
for large $|s|$. This together with (\ref{2.7}) yields that
\begin{equation*}
\frac{1}{2\pi}\int_0^{2\pi}|\log^+|f_1(Re^{i\phi})||d\phi\le \epsilon R^2
\end{equation*}
for large $R.$ We then deduce from (\ref{2.6}) that with $s=re^{i\theta}$,
\begin{eqnarray}\label{estimate of f_1}
& &\frac{1}{2\pi}\int_{E}\log|\frac{f_1(s+\eta+\tau)}{f_1(s+\eta)}|d\theta\\ \nonumber
& &\le \frac{2\tau\epsilon R^3}{(R-r-d\tau)^2}+\sum_{k=1}^{M}\frac{1}{2\pi}\int_E\log|\frac{R^2-\bar{a_k}(s+\eta)}{R(s+\eta-a_k)}|d\theta.
\end{eqnarray}
To estimate the summation in the above inequality, we notice by (\ref{2.5}) that $|g(s)|\ge 1$ in $|s|\le R$ and thus $g$ has no zeros in $|s|\le r<R$. Also, $g$ has only one pole $a_k-\eta$ in $|s|\le r$ if $|a_k-\eta|\le r$ and has no poles in $|s|\le r$ if $|a_k-\eta|>r$. Hence, applying the Jensen formula (\ref{jensen}) to the function $g(s):=\frac{R^2-\bar{a_k}(s+\eta)}{R(s+\eta-a_k)}$ in $|s|\le r$, we obtain that
\begin{eqnarray*}
& &\log|g(0)|=\log|\frac{R^2-\bar{a_k}\eta}{R(\eta-a_k)}|\\
& &=\frac{1}{2\pi}\int_0^{2\pi}\log|\frac{R^2-\bar{a_k}(re^{i\theta}+\eta)}{R(re^{i\theta}+\eta-a_k)}|d\theta+\log^+\frac{r}{|a_k-\eta|}
\end{eqnarray*}
and thus, noting (\ref{2.5}) again, that with $s=re^{i\theta}$,
\begin{eqnarray*}
& &\sum_{k=1}^{M}\frac{1}{2\pi}\int_E\log|\frac{R^2-\bar{a_k}(s+\eta)}{R(s+\eta-a_k)}|d\theta\\
& &\le \sum_{k=1}^{M}\frac{1}{2\pi}\int_0^{2\pi}\log|\frac{R^2-\bar{a_k}(s+\eta)}{R(s+\eta-a_k)}|d\theta\\
& &=\sum_{k=1}^{M}\log|\frac{R^2-\bar{a_k}\eta}{R(\eta-a_k)}|-\sum_{k=1}^{M}\log^+\frac{r}{|a_k-\eta|}\\
& &\le \sum_{k=1}^{M}\log \frac{R+d\tau}{|\eta-a_k|}.
\end{eqnarray*}
Since the zeros of $f_1(s)$ are $a_k$ ($k=1, 2,\cdots$),  the zeros of $f_1(s+\eta)$ are $a_k-\eta$ ($k=1, 2,\cdots$). Note also that when $|a_k|\le R$, $|a_k-\eta|\le R+|\eta|\le R+d\tau$ by the definition of $\eta$. Thus, we deduce, in view of (\ref{sum}) and the fact that $f_1$ does not vanish in $|s|<2d\tau$, that
\begin{eqnarray*}
& &\sum_{k=1}^{M}\log \frac{R+d\tau}{|\eta-a_k|}\\
& &\le\sum_{|a_k-\eta|\le R+d\tau}\log \frac{R+d\tau}{|a_k-\eta|}\\
& &=\int_0^{R+d\tau}\frac{n(t, 0; f_1(s+\eta))-n(0, 0; f_1(s+\eta))}{t}dt\\
& &\le \int_0^{R+d\tau}\frac{n(t+d\tau, 0; f_1(s))}{t}dt\\
& &=\int_{d\tau}^{R+2d\tau}\frac{n(t, 0; f_1(s))}{t-d\tau}dt \\
& &=\int_{2d\tau}^{R+2d\tau}\frac{n(t, 0; f_1(s))}{t}(1+\frac{d\tau}{t-d\tau})dt\\
& &\le 2N(R+2d\tau, 0; f).
\end{eqnarray*}
This together with (\ref{estimate of f_1}) yields that with $s=re^{i\theta}$,
\begin{eqnarray*}
& &\frac{1}{2\pi}\int_{E}\log|\frac{f_1(s+\eta+\tau)}{f_1(s+\eta)}|d\theta\\
& &\le \frac{2\tau\epsilon R^3}{(R-r-d\tau)^2}+2N(R+2d\tau, 0; f).
\end{eqnarray*}

This same estimate holds for all the integrals in (\ref{2.3}) by the same argument. Thus, in view of  (\ref{integrated counting}), (\ref{2.3}) and (\ref{2.4}) we obtain that for a small number $\epsilon>0$ and a constant $C>0$,
\begin{eqnarray}\label{estimate of N}
& &cr\le N(r, a, F) \\
& &\le \frac{\epsilon R^3}{(R-r-d\tau)^2}+C\bigl(N(R+2d\tau, 0; f)+N(R+2d\tau, \infty; f)+\log r \bigr)    \nonumber
\end{eqnarray}
for $r>e$ and large $R>r+d\tau$, where $c>0$ is the constant in (\ref{integrated counting}).

Next, we bound the right hand side of the estimate (\ref{estimate of N}) in terms of $n(r, 0; f)+n(r, \infty; f)$ so that we can make use of the assumption (\ref{2.1}) to arrive at the desired contradiction mentioned in the beginning of the proof.  (It would be tempting to try to use the definition of the integrated counting function to bound $N(r, 0; f)+N(r, \infty; f)$ using $n(r, 0; f)+n(r, \infty; f)$ directly; but it would only yield the rough estimate that
$
N(r, 0; f)+N(r, \infty; f)\le \bigl(n(r, 0; f)+n(r, \infty; f)\bigr)\log r+O(\log r),
$
which does not serve our purpose.) To this end, fix a large positive number $R_1.$ Let $a_1, a_2,\cdots,a_p$ be the zeros and  $b_1, b_2,\cdots,b_q$ be the poles of $f$ in $|s|\le R_1$. Thus, $p=n(R_1, 0; f)$ and $q=n(R_1, \infty, f)$. Applying the Cartan theorem (see (\ref{cartan})) to the zeros and poles of $f$, respectively,  with $h=\frac{1}{32}R_1$, we obtain a collection ${\cal D}$ of the disks with the sum of the radii equal to $4h=\frac{1}{8}R_1$ such that for all $z\not\in {\cal D}$ we have
\begin{equation}\label{disk1}
\Pi_{k=1}^p|s-a_k|>(\frac{h}{e})^p
\end{equation}
and
\begin{equation*}
\Pi_{k=1}^q|s-b_k|>(\frac{h}{e})^q.
\end{equation*}
Since the sum of the radii of the disks in ${\cal D}$ is $\frac{1}{8}R_1$, it is easy to check that those disks in ${\cal D}$ cannot fill up the annulus $\cal{G}$: $\frac{1}{16}R_1\le |s|\le \frac{1}{4}R_1$. Therefore, there must be a point $s_0\in \cal{G}$ outside all the above disks, i.e., $s_0\not\in {\cal D}$ with $f(s_0)\not=0, \infty$. Note that the estimate (\ref{estimate of N}) was obtained for any $f$ that satisfies the conditions of Theorem 1.2 and the assumption (\ref{con3}), which are satisfied by $f(s+s_0)$ as well. Therefore, applying the proved estimate (\ref{estimate of N}) to $f(s+s_0)$, we obtain that for large $R$ and $r$ with $R>r+d\tau$,
\begin{eqnarray}\label{estimate1 of N}
& &cr\le N(r, a, F(s+s_0))  \nonumber\\
& &\le \frac{\epsilon R^3}{(R-r-d\tau)^2}+C\bigl(N(R+2d\tau, 0; f(s+s_0))+ \nonumber\\
& &+N(R+2d\tau, \infty; f(s+s_0))+\log r\bigr).
\end{eqnarray}
Note that $f(s)$ vanishes at $a_k$ if and only if $f(s_0+s)$ vanishes at $a_k-s_0$. Denote $\rho=\frac{3}{4}R_1$, we have, in view of (\ref{sum}), that
\begin{eqnarray*}
& &N(\rho, 0; f(s_0+s))\\
& &=\sum_{|a_k-s_0|\le \rho}\log\frac{\rho}{|a_k-s_0|}\\
& &=\log \frac{\rho^{n_1}}{\Pi_{|a_k-s_0|\le \rho}|a_k-s_0|}\\
& &=\log \frac{\rho^{n_1}\Pi_{|a_k|\le R_1, |a_k-s_0|>\rho}|a_k-s_0|}{\Pi_{|a_k|\le R_1}|a_k-s_0|},
\end{eqnarray*}
where $n_1$ is the number of the points $a_k$ satisfying $|a_k-s_0|\le \rho$; these $a_k$'s satisfy that $|a_k|\le \rho+|s_0|\le \frac{3}{4}R_1+\frac{1}{4}R_1=R_1$. Let $n_2$ be the number of $a_k$ satisfying that $|a_k-s_0|>\rho$ and $|a_k|\le R_1$. It is then clear that $n_1+n_2=p=n(R_1, 0; f).$ Thus, by virtue of (\ref{disk1}) and the fact that
$|a_k-s_0|\le 2R_1$ for any $k\le p$, we deduce that
\begin{eqnarray*}
& &N(\rho, 0; f(s_0+s))\\
& &\le \log \frac{\rho^{n_1}(2R_1)^{n_2}}{(\frac{R_1}{32e})^p}\\
& &\le \log \frac{(2R_1)^p}{(\frac{R_1}{32e})^p}=n(R_1, 0; f)\log (64e).
\end{eqnarray*}
In the exactly same way, we also have that
\begin{equation*}
N(\rho, \infty; f(s_0+s))\le n(R_1, \infty; f)\log (64e).
\end{equation*}
Thus, we have that
\begin{eqnarray}\label{counting4}
& &N(\frac{3}{4}R_1, 0; f(s_0+s))+N(\frac{3}{4}R_1, \infty; f(s_0+s)) \nonumber\\
& &\le \bigl(n(\frac{3}{4}R_1, 0; f)+n(\frac{3}{4}R_1, \infty; f)\bigr)\log (64e).
\end{eqnarray}
Fix a large $n$. Take $R_1=\frac{4}{3}r_n,$ where $\{r_n\}$, which tends to $+\infty$, is the sequence in (\ref{2.1}). Then, take $r=\frac{1}{2}r_n$ and $R=\frac{3}{4}R_1-2d\tau$ in (\ref{estimate1 of N}). Combining (\ref{estimate1 of N}) and (\ref{counting4}) and (\ref{estimate1 of N}) yields that
\begin{eqnarray*}
& &cr_n\le\frac{\epsilon (r_n-2d\tau)^3}{(\frac{1}{2}r_n-3d\tau)^2}+C\bigl((n(r_n, 0; f)+n(r_n, \infty; f))\log (64e)+\log r_n\bigr)\\
& &\le\frac{\epsilon (r_n-2d\tau)^3}{(\frac{1}{2}r_n-3d\tau)^2}+C\bigl(\frac{1}{n}r_n\log (64e)+\log r_n\bigr),\\
\end{eqnarray*}
where $c, C>0$ are fixed constants, $\epsilon>0$ can be made arbitrarily small and $n$ can be made arbitrarily large, which is clearly impossible. This finally completes the proof of the theorem.
\qed

\bigskip
\noindent
\textit{Department of Mathematics and Statistics\\
Florida International University\\
Miami, FL 33199 USA\\
libaoqin@fiu.edu}

\end{document}